\documentclass{my-elsarticle}

\journal{Alexandria Engineering Journal}

\usepackage{amsthm,amsmath,amssymb}
\usepackage{graphicx}
\usepackage{url}
\usepackage{a4wide}

\newtheorem{theorem}{Theorem}

\newtheorem{definition}[theorem]{Definition}

\begin{document}

\begin{frontmatter}

% --------------------------------------

\title{Control of COVID-19 dynamics through a fractional-order model}

% --------------------------------------

\author[aff4]{Samia Bushnaq}
\ead{S.Bushnaq@psut.edu.jo}
\ead[url]{https://orcid.org/0000-0002-2427-7704}

\author[aff3]{Tareq Saeed}
\ead{tsalmalki@kau.edu.sa}
\ead[url]{https://orcid.org/0000-0002-0170-5286}

\author[aff2]{Delfim F. M. Torres\corref{cor1}}
\ead{delfim@ua.pt}
\ead[url]{https://orcid.org/0000-0001-8641-2505}
\cortext[cor1]{Corresponding author}

\author[aff1]{Anwar Zeb}
\ead{anwar@cuiatd.edu.pk}
\ead[url]{https://orcid.org/0000-0002-5460-3718}

% --------------------------------------

\address[aff4]{Department of Basic Sciences,
Princess Sumaya University for Technology,
11941 Amman, Jordan}

\address[aff3]{Department of Mathematics,
King Abdulaziz University,
41206 Jeddah, Kingdom Saudi Arabia}

\address[aff2]{Center for Research and Development in Mathematics and Applications (CIDMA),\\
Department of Mathematics, University of Aveiro, 3810-193 Aveiro, Portugal}

\address[aff1]{Department of Mathematics,
COMSATS University Islamabad,
22060 Abbottabad, Pakistan}

% --------------------------------------

\begin{abstract}
We investigate, through a fractional mathematical model,
the effects of physical distance on the SARS-CoV-2 virus
transmission. Two controls are considered in our model
for eradication of the spread of COVID-19: media education,
through campaigns explaining the importance  of social distancing,
use of face masks, etc., towards all population, while the second
one is quarantine social isolation of the exposed individuals.
A general fractional order optimal control problem, and
associated optimality conditions of Pontryagin type,
are discussed, with the goal to minimize the number of susceptible
and infected while maximizing the number of recovered. The extremals
are then numerically obtained.
\end{abstract}

% --------------------------------------


\end{frontmatter}

% --------------------------------------

\section{Introduction}

The availability of easy-to-use precise estimation models are
essential to get an insight into the effects of transferable
infectious diseases. In outbreak diseases, policy makers and
institutions make decisions based on forecasting models to decide on
future policies and to check the efficiency of existing policies
\cite{c1}.

Coronaviruses are a group of viruses that can be
transmitted between humans, livestock and wild animals.
Person to person spread of COVID-19 happens through
close contact, up to six feet.
This group of viruses mainly affects the hepatic, neurological
and respiratory systems \cite{c2,c3,c4}.

In the end of 2019, the World Health Organization (WHO)
reported a novel coronavirus in China,
which causes severe damage to the respiratory system. The virus
was first found in Wuhan city, and was named as severe acute
respiratory syndrome coronavirus 2 (SARS-CoV-2) \cite{c5}.
After the outbreak of the virus, the Chinese government
put several cities on lock down \cite{c6,c7}. However,
the number of affected people increased daily within China
and in other countries. In March 11, 2020, COVID-19 was declared
as a global pandemic by WHO. At the time we write these lines,
January 6, 2021, approximately 87 million people were infected
with 1.9 million of deaths, worldwide \cite{worldometer}.

Recently, many mathematical models were proposed to understand
disease transmission and project handful controls: see, e.g.,
\cite{LP,MR4093642,MR4173153} and references therein.
Simultaneously, all health organizations are trying to drive
the most lethal infectious diseases towards eradication, using
educational and enlightenment campaigns, vaccination, treatment, etc.
However, many of these infectious diseases will become eventually
endemic because of interventions to mitigate the spread in time
and lack of adequate policies. For control of infectious diseases,
proactive steps are required, specially for diseases having vaccine and cure.
Indeed, some times it is more difficult to control the spread of an infectious
disease than to cure it. Regarding COVID-19, several vaccines begin to be available
\cite{tracker}.

Optimal control theory is a branch of mathematical optimization that
deals with finding a control for a dynamical system, over a period
of time, such that an objective function is minimized or maximized.
Along the years, optimal control theory has found applications in
several fields, containing process control, aerospace, robotics,
economics, bio-engineering, management sciences, finance, and
medicine \cite{MR3999699,MR4110649,MR3988048}. In particular, the
study of epidemic models is strongly related to the study of control
strategies, as screening and educational campaigns \cite{Castillo},
vaccination \cite{Brandeau}, and resource allocation \cite{Becker}.

In the current pandemic situation of COVID-19, due to best
presentation of memory effects and its usefulness in many different
and widespread phenomena
\cite{Anastasio,Agrawal,Diethelm,Diethelm1,Agrawal1,Tricaud,Biswas},
fractional (non-integer order) models are receiving the attention of
many researchers: see, e.g.,
\cite{r1,MyID:473,r2,r3,r4,r5,r6,r7,r8,r9,r10,r11,r12}. Here, with
the purpose to control the current pandemic, we follow two control
variables, in the form of media (education) campaigns, social
distance, and use of mask for protection of susceptible individuals;
and quarantine (social isolation) for the exposed. For a general
non-integer order optimal control problem, necessary optimality
conditions are presented, with the help of Caputo derivatives.
One of the great advantages of the
Caputo fractional derivative is that it allows traditional initial
or boundary conditions to be included in the formulation of the
problem. More concretely, we minimize the number of susceptible and
infected, while maximizing the number of recovered 
population from COVID-19. The optimal levels of the proposed two
controls are characterized using the fractional version of
Pontryagin's maximum principle. The resulting optimality system is
then solved numerically with \textsf{Matlab}.

The rest of the paper is arranged as follows. In Section~\ref{sec2}, we
present our fractional mathematical model. Section~\ref{sec3}
recalls the fundamental definitions and the main result of fractional optimal control.
We then derive an optimal control problem in Section~\ref{sec4},
while parameter estimation with numerical results are discussed
in Section~\ref{sec5}. We finish with some remarks
and conclusions, in Section~\ref{sec6}.

% --------------------------------------

\section{The fractional model}
\label{sec2}

Our model consists of four classes: $S(t)$, which represents the
vulnerable individuals (healthy people but who may get the disease
in a near future); $E(t)$, representing the exposed population or
individuals who are infected but not yet infectious; the group
$I(t)$, devoted to the population who are confirmed infected
(individuals who have contracted the disease and are now sick with
it and are infectious); the group $R(t)$, defined as the recovered
population (individuals who have recovered from COVID-19).
For the dynamics of this base model, see \cite{r2}. 
Thus, the fractional order model we consider here is given by
\begin{equation}
\label{eq1.1}
\left\{
\begin{array}{rcl}
_{0}^{C}D_{t}^{\alpha}S(t) &=&\Lambda-\beta_{1}S(t)E(t)
-\beta_{2}S(t)I(t)-\mu S(t)+\tau R(t),\cr _{0}^{C}D_{t}^{\alpha}E(t)
&=&\beta_{1}S(t)E(t)+\beta_{2}S(t)I(t) -(\mu+\rho)E(t),\cr
_{0}^{C}D_{t}^{\alpha}I(t) &=&\rho E(t)-(\gamma +d+\mu)I(t),\cr
_{0}^{C}D_{t}^{\alpha}R(t) &=&\gamma I(t)-(\mu+\tau)R(t),
\end{array}
\right.
\end{equation}
where $\Lambda$ is the recruitment
rate, $\beta_{1}$ and $\beta_{2}$ are the incidence rates, $\tau$ is
the relapse rate, $\mu$ is the natural death rate, $\rho$ the rate
at which the exposed population of COVID-19 join the infectious
class, $\gamma$ the recovery rate of infected population, and $d$ is
the death rate of infected class due to the SARS-CoV-2 virus. The
total population $N(t)$ is given, at each instant of time, by
\begin{equation}
\label{eq1.2}
N(t)=S(t)+E(t)+I(t)+R(t).
\end{equation}
By adding all the equations of system (\ref{eq1.1}), we have
\begin{equation}
\label{eq1.3} _{0}^{C}D_{t}^{\alpha}N(t) =\Lambda-\mu
N(t)-dI(t)\leq\Lambda-\mu N(t).
\end{equation}

% --------------------------------------

\section{Basics of fractional control theory}
\label{sec3}

In this section we recall the basic definitions of Caputo
fractional calculus and the central result of
fractional optimal control theory \cite{Agrawal1,MR3787674},
which are required for the coming sections.

\begin{definition}
For $f \in C^m$, $m \in \mathbb{N}$, the left-sided Caputo
fractional derivative is given by
\begin{equation}
\label{eq2.1}
{_{a}^{C}D_{t}^{\alpha}}f(t)=\frac{1}{\Gamma(m-\alpha)}\int_{a}^{t}
(t-\tau)^{m-\alpha-1}\left(\frac{d}{d\tau}\right)^{m}f(\tau)d\tau,
\end{equation}
while the right-sided Caputo fractional derivative is given by
\begin{equation}
\label{eq2.2}
{_{t}^{C}D_{b}^{\alpha}}f(t)=\frac{1}{\Gamma(m-\alpha)}{\int_{t}^{b}}
(\tau-t)^{m-\alpha-1}\left(-\frac{d}{d\tau}\right)^{m}f(\tau)d\tau,
\end{equation}
where  $\alpha$ stands for order of the derivative, $m-1<\alpha\leq
m$.
\end{definition}

\begin{definition}
For $f$ an integrable function, the left-sided
Riemann--Liouville fractional derivative is defined by
\begin{equation}
\label{eq2.3}
{_{a}D_{t}}^{\alpha}f(t)=\frac{1}{\Gamma(m-\alpha)}\left( \frac{d}
{dt}\right) ^{m}\int_{a}^{t}(t-\tau)^{m-\alpha-1}f(\tau)d\tau,
\end{equation}
while the right-sided Riemann--Liouville derivative of $f$ is given by
\begin{equation}
\label{eq2.4}
{_{t}D_{b}^{\alpha}}f(t)=\frac{1}{\Gamma(m-\alpha)}\left( -\frac{d}
{dt}\right) ^{m}{\int_{t}^{b}}(\tau-t)^{m-\alpha-1}f(\tau)d\tau,
\end{equation}
where  $\alpha$ is the order of the derivative with $m-1<\alpha \leq m$,
$m \in \mathbb{N}$.
\end{definition}

Our control system is described by a fractional differential system (FDS)
with a given/fixed initial condition as follows:
\begin{equation}
\label{eq2.5}
\begin{cases}
_{0}^{C}D_{t}^{\alpha}X(t)=f\left(  X(t),u(t),t\right), \\
X(0)=X_{0},
\end{cases}
\end{equation}
where $\alpha\in(0,1]$, the $n$-dimensional $X(t)$ is the state vector,
$f$ is a given vector-valued function, $t\in\lbrack0,t_{f}]$ with $t_{f}>0$
the ending time of the control process, and $m$-dimensional $u(t)$
is the control vector. A fractional optimal control problem consists
to minimize or maximize a performance index
\begin{equation}
\label{eq2.6}
J\left[ u(\cdot)\right]  =\theta\left(  X(t_{f}),t_{f}\right)
+ \int\limits_{0}^{t_{f}} \phi\left( X(t),u(t),t\right) dt
\end{equation}
subject to the control system \eqref{eq2.5} 
(see, e.g., \cite{MR3787674,MR3443073}). 
Functions $\theta$ and $\phi$ will be specified in Section~\ref{sec4}.
Note that here $t_{f}$ is fixed but $X(t_{f})$ is free. 
For finding the optimal control law
$u(t)$ solution to the optimal control problem
\eqref{eq2.5}--\eqref{eq2.6}, we use the fractional version of
Pontryagin maximum principle, which coincides with the classical
Pontryagin maximum principle when $\alpha=1$:

\begin{theorem}[See, e.g., \cite{MR3787674,MR3443073}]
\label{PMP} For the optimality of \eqref{eq2.5}--\eqref{eq2.6}, a
necessary condition is given by
\begin{equation*}
\begin{cases}
\frac{\partial\phi}{\partial u}\left( X(t),u(t),t\right)
+\lambda^{T} \frac{\partial f}{\partial u}\left(  X(t),u(t),t\right) = 0, \\
_{0}^{C}D_{t}^{\alpha}X(t)=f\left(  X(t),u(t),\left(  t\right)  \right), \quad X(0)=X_{0},\\
_{t}D_{t_{f}}^{\alpha}\lambda(t)=\frac{\partial\phi}{\partial X}\left( X(t),u(t),t\right)
+\lambda ^{T}\frac{\partial f}{\partial X}\left(  X(t),u(t),t\right),
\quad \lambda(t_{f})=\frac{\partial\theta}{\partial X}\left(X(t_{f}),t_{f}\right).
\end{cases}
\end{equation*}
\end{theorem}

In the next section, we compute the optimal control strategy
for the fractional order COVID-19 model, which is
a hot topic in current times.

% --------------------------------------

\section{Fractional-order model with controls}
\label{sec4}

We implement an optimal control technique to the fractional order
model (\ref{eq1.1}). With the purpose to control the spread of the
COVID-19 pandemic in the world, we use two control variables in the
form of media (educational) campaigns, social distance, and use of
masks --- the control $u_{1}(t)$ --- applied to the susceptible
class; and quarantine (social isolation) --- the control $u_{2}(t)$
--- applied to the exposed class. 
Then, the new system with controls is given by
\begin{equation}
\label{eq3.1}
\begin{cases}
_{0}^{C}D_{t}^{\alpha}S(t) =\Lambda-\beta_{1}S(t)E(t)
-\beta_{2}S(t)I(t)-\mu S(t)+\tau R(t)-u_{1}(t) S(t),\cr
_{0}^{C}D_{t}^{\alpha}E(t) =\beta_{1}S(t)E(t)+\beta_{2}S(t)I(t)
-(\mu+\rho)E(t)-u_{2}(t)E(t),\cr _{0}^{C}D_{t}^{\alpha}I(t) =\rho
E(t)-(\gamma+d+\mu)I(t)+(1-p)u_{2}(t)E(t),\cr
_{0}^{C}D_{t}^{\alpha}R(t) =\gamma I(t)-(\mu+\tau)R(t)+u_{1}(t)
S(t)+pu_{2}(t)E(t),
\end{cases}
\end{equation}
where the fractional order $\alpha$ is a real number in the interval $(0,1]$
and $p$ can be interpreted as the probability
of infected individuals to recover by quarantine. 
In vector form, the system (\ref{eq3.1}) can be written as
\begin{equation}
\label{eq3.2}
_{0}^{C}D_{t}^{\alpha}{X(t)}=f\left( X(t), u(t)\right),
\end{equation}
where $X(t)=\left( S(t), E(t), I(t), R(t)\right)$ represents the
state-vector and $u(t)=\left( u_{1}(t), u_{2}(t)\right)$ stands
for the control-vector.

Our optimal control problem consists to minimize the spread of
COVID-19 and maximize the number of recovered population.
The following objective functional is defined with this purpose:
\begin{equation}
\label{eq3.3}
J[u(\cdot)]=
A_{3} S(t_f)+A_{4} E(t_f)+
\int_{0}^{t_{f}}
A_{1}I(t)-A_{2}R(t)+\frac{1}{2}\left(r_{1}u_{1}^{2}(t)+r_{2}u_{2}^{2}(t)\right) dt
\longrightarrow \min,
\end{equation}
where the positive weights $A_{i}$, $i=1,2,3,4$, and
$r_{i}$, $i=1,2$, are used to balance the control factors.
The objective functional \eqref{eq3.3} is a particular case
of the general form \eqref{eq2.6} discussed in Section~\ref{sec3},
and can be written as
\begin{equation}
\label{eq3.4}
J[u(\cdot)]= \theta\left(X(t_{f})\right)
+\int_{0}^{t_{f}}\phi\left(X(t),u(t)\right)dt
\end{equation}
with $\theta\left(X(t_{f})\right)= A_{3}S(t_f)+A_{4}E(t_f)$ and
\begin{equation*}
\phi\left(X(t),u(t)\right)
=A_{1}I(t)-A_{2}R(t)+\frac{1}{2}\left(r_{1}u_{1}^{2}(t)+r_{2}u_{2}^{2}(t)\right).
\end{equation*}
Similar functionals \eqref{eq3.4}
to be optimized, e.g. for optimal control problems in the combat of Zika and Ebola, 
have been previously considered in the literature,
see \cite{MyID:392,MyID:364} and references therein.
By using Theorem~\ref{PMP}, the following necessary optimality
conditions can be written: the control system and its initial
condition,
\begin{equation}
\label{eq3.5}
\begin{cases}
_{0}^{C}D_{t}^{\alpha}X=f\left(  X,u\right), \\
X(0)=X_{0},
\end{cases}
\end{equation}
the adjoint system and its transversality condition,
\begin{equation}
\label{as:tc}
\begin{cases}
_{t}D_{t_{f}}^{\alpha}\lambda(t)=\frac{\partial\phi}{\partial X}
+\lambda ^{T}\frac{\partial f}{\partial X}, \\
\lambda(t_{f})=\left.\frac{\partial\theta}{\partial X}\right|_{t_{f}},
\end{cases}
\end{equation}
and the stationary condition
\begin{equation}
\frac{\partial\phi}{\partial u}
+\lambda^{T} \frac{\partial f}{\partial u} = 0,
\end{equation}
where $\lambda(t)=\left(\lambda_{1}(t),\lambda_{2}(t),\lambda_{3}(t),\lambda_{4}(t)\right)$
and $f=\left(f_{1},f_{2},f_{3},f_{4}\right)$ with
\begin{equation*}
\begin{split}
f_1 &= \alpha-\beta_{1}S(t)E(t)-\beta_{2}S(t)I(t)-\mu S(t)+\tau R(t)-u_{1}(t) S(t),\\
f_2 &= \beta_{1}S(t)E(t)+\beta_{2}S(t)I(t)-(\mu+\rho)E(t)-u_{2}(t)E(t),\\
f_3 &= \rho E(t)-(\gamma+d+\mu)I(t)+(1-p)u_{2}(t)E(t),\\
f_4 &= \gamma I(t)-(\mu+\tau)R(t)+u_{1}(t) S(t)+pu_{2}(t)E(t).
\end{split}
\end{equation*}
The adjoint system of Pontryagin's maximum principle asserts that
\begin{equation}
\label{eq3.6}
\begin{cases}
_{t}D_{t_f}^{\alpha}\lambda_1(t)
=-\lambda_1\beta_{1}E(t)-\lambda_1\beta_{2}I(t)
-\lambda_1\mu+\lambda_2\beta_{1}E(t)+\lambda_2\beta_{2}I(t),\cr
_{t}D_{t_f}^{\alpha}\lambda_2(t)
=-\lambda_1\beta_{1}S(t)+\lambda_2\beta_{1}S(t)
-\lambda_2\mu-\lambda_2\rho+\lambda_3\rho,\cr
_{t}D_{t_f}^{\alpha}\lambda_3(t)
=A_1-\beta_{2}\lambda_1S(t)+\beta_{2}\lambda_2S(t)
-(\gamma+d+\mu)\lambda_3+\gamma\lambda_4,\cr
_{t}D_{t_f}^{\alpha}\lambda_4(t)
=\tau\lambda_1-\lambda_4(\tau+\mu)-A_2,
\end{cases}
\end{equation}
subject to the transversality conditions
\begin{equation}
\label{eq3.7}
\begin{cases}
\lambda_1(t_f) = A_3,\cr
\lambda_2(t_f) = A_4,\cr
\lambda_3(t_f) = 0,\cr
\lambda_4(t_f) = 0.
\end{cases}
\end{equation}
The optimal control variables are given by the stationary conditions:
\begin{equation}
\label{eq3.8}
\begin{cases}
u_{1}(t)=\frac{(\lambda_1(t)- \lambda_4(t))S(t)}{r_{1}},\cr
u_2(t)=\frac{(\lambda_2-(1-\rho)\lambda_3-\rho\lambda_4)E(t)}{r_{2}}.
\end{cases}
\end{equation}
These analytic necessary optimality 
conditions are solved numerically in Section~\ref{sec5}.

% --------------------------------------

\section{Numerical simulations}
\label{sec5}

To illustrate the theoretical results presented in previous
sections, here we use numerical simulations. For this purpose, a
program was developed in \textsf{Matlab} to integrate the necessary
optimality conditions and, with the help of a number of simulations,
a detailed output is comprehensively verified. As explained in
Section~\ref{sec4}, we obtain the optimality system for the proposed
optimal control problem from the state and adjoint equations subject
to suitable boundary conditions: the initial conditions 
$X(0)=X_{0}$ on the state variables, see \eqref{eq3.5};
and the terminal conditions on the adjoint variables provided 
by the transversality conditions, see \eqref{eq3.7}. Furthermore, 
we obtain the optimal control strategies from the stationary system, 
see \eqref{eq3.8}. We use a forward time/backward 
space finite-difference numerical method. Beginning
with an initial guess for the adjoint variables, a forward time and
backward space finite-difference method is used to solve the state
equations. The key is to rewrite the control system \eqref{eq3.2}
into the equivalent integral form
\[
X(t)=X(0)+\frac{1}{\Gamma(\alpha)}
{\displaystyle\int\limits_{0}^{t}}
(t-\tau)^{\alpha-1}f(X(\tau),u(\tau))d\tau
\]
and then use the generalized Adams-type predictor-corrector method
\cite{Diethelm,Diethelm1} for solution. Further, these state values
are used for the solution of the adjoint equations by a backward
time and forward space finite-difference method, because of the
transversality conditions. System \eqref{as:tc} is written, in an
equivalent way, as the integral equation
\[
\lambda(t)= \left.\frac{\partial\theta}{\partial X}\right|_{t_{f}}
+\frac{1}{\Gamma(\alpha)} {\displaystyle\int\limits_{t}^{t_{f}}}
(\tau-t)^{\alpha-1}\left[  \frac{\partial\phi}{\partial
X}+\lambda^{T} \frac{\partial f}{\partial X}\right]  d\tau.
\]
Using a steepest-method to generate a successive approximation of
the optimal control form, we continue iterating until convergence is
achieved. For illustrative purposes, take the initial values as
$S(0)=220$, $E(0)=100$, $I(0)=3$, $R(0)=0$ and parameter values as
$\Lambda=0.271$, $\beta_1=0.00035$, $\beta_2=0.00040$, $\mu=0.001$,
$\rho=0.00580$, $\gamma=0.007$, $\tau=0.002$, $p=0.3$, and $d=0.00025$.

\begin{figure}[hbtp]
\begin{center}
\includegraphics[scale=0.70]{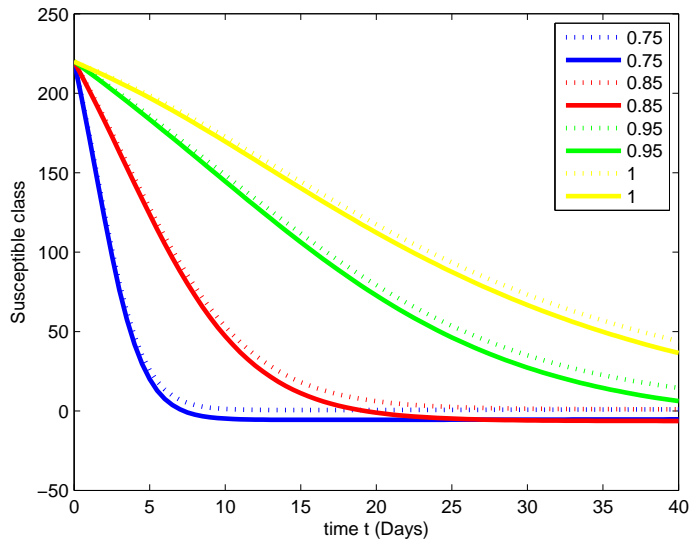}
\caption[Figure]{The susceptible population $S(t)$, with and without
controls, respectively solid and doted lines,
for $\alpha=0.75,0.85,0.95,1$.} \label{Fig.1}
\end{center}
\end{figure}

In Figure~\ref{Fig.1}, we plot the susceptible population of systems
(\ref{eq1.1}) and (\ref{eq3.1}). The doted lines denote the
population of class $S$ in the uncontrolled system (\ref{eq1.1}),
without controls, while the solid lines denote the population of
$S(t)$ in the controlled system (\ref{eq3.1}), under optimal
controls for $\alpha=0.75$, $0.85$, $0.95$ and $1$.
% ----------------------------
% ----------------------------

\begin{figure}[hbtp]
\begin{center}
\includegraphics[scale=0.70]{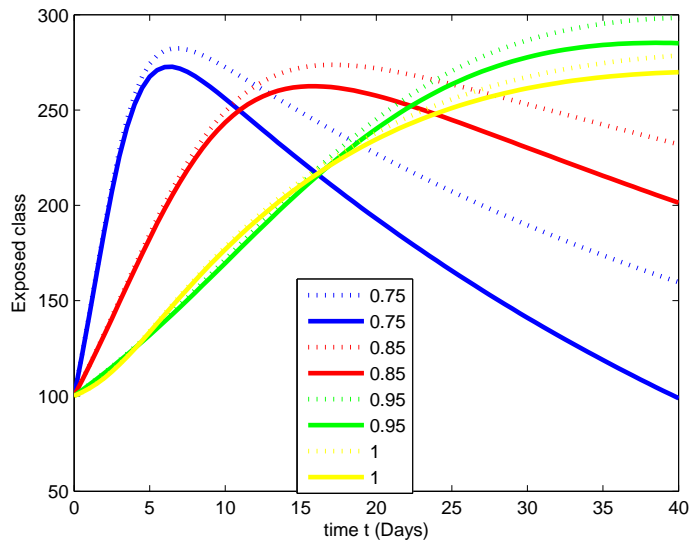}
\caption[Figure]{The exposed class $E(t)$ of individuals, with and
without controls, respectively solid and doted lines,
for $\alpha=0.75,0.85,0.95,1$.} \label{Fig.2}
\end{center}
\end{figure}

Figure~\ref{Fig.2} represents the exposed population of both systems
(\ref{eq1.1}) and (\ref{eq3.1}). The doted lines show that there are
more exposed individuals when no control measures are implemented.
% ----------------------------
% ----------------------------

\begin{figure}[hbtp]
\begin{center}
\includegraphics[scale=0.70]{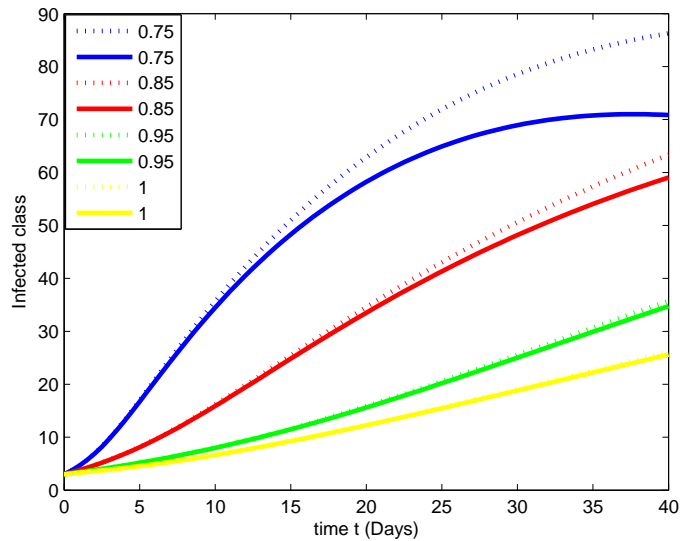}
\caption[Figure]{Infected population for systems with and without
controls, respectively solid and doted lines,
for $\alpha=0.75,0.85,0.95,1$.} \label{Fig.3}
\end{center}
\end{figure}

Figure~\ref{Fig.3} illustrates the infectious population $I(t)$ of
system (\ref{eq1.1}), without any control, and that of system
(\ref{eq3.1}) with controls. The doted lines make it clear that
there are more infectious individuals when no control is
implemented.
% ----------------------------
% ----------------------------

\begin{figure}
\begin{center}
\includegraphics[scale=0.70]{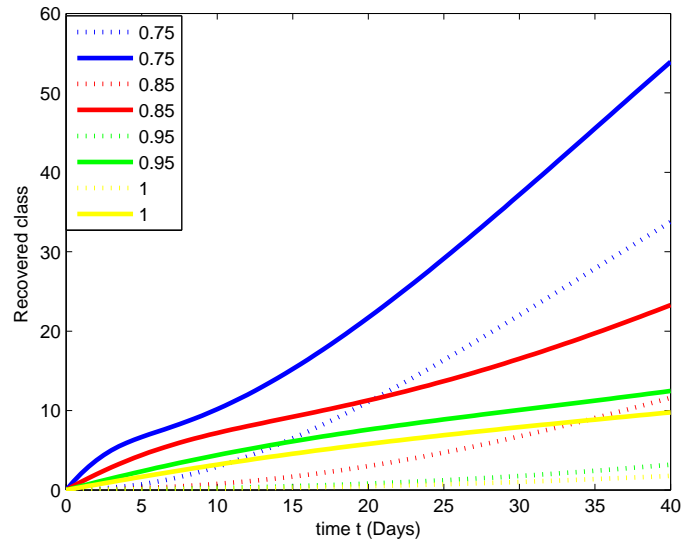}
\caption[Figure]{Recovered individuals for systems with and without
controls, respectively solid and doted lines,
for $\alpha=0.75,0.85,0.95,1$.} \label{Fig.4}
\end{center}
\end{figure}

Finally, Figure~\ref{Fig.4} illustrates the recovered population $R(t)$.
We see that there are more recovered individuals
in the case one uses optimal control theory
(because there is less susceptible, exposed and infected).

% ----------------------------

% --------------------------------------

\section{Conclusion}
\label{sec6}

The current pandemic situation due to COVID-19
affects the whole world on an unprecedented scale.
In this work, we implemented optimal control techniques to the
COVID-19 pandemic through a fractional order model. For the
eradication of virus spread throughout the world, we applied two
controls in the form of media (education) campaigns, social
distance, use of masks and protection for the susceptible class; and
quarantine (social isolation) for the exposed individuals. We
discussed necessary optimality conditions for a general fractional
optimal control problem, whose fractional system is described in the
Caputo sense while the adjoint system involves Riemann--Liouville
derivatives. In the COVID-19 setting, we minimize the number of
susceptible and infected population, while maximizing the number of
recovered population from SARS-CoV-2 virus. Using the
fractional version of Pontryagin's maximum principle, we
characterize the optimal levels of the proposed controls.
The resulting optimality system is solved numerically
in the \textsf{Matlab} numerical computing environment.
Our numerical experiments were based on data of \cite{MyID:473}. 
In a future work, we plan to use real data of Africa, USA and UK.

% --------------------------------------

\section*{Acknowledgement}

D.F.M.T. was supported by The Center for Research and Development in
Mathematics and Applications (CIDMA) through FCT, project
UIDB/04106/2020. The authors are very grateful
to two Reviewers for several useful comments, suggestions and questions,
which helped them to improve the original manuscript.

% --------------------------------------

% --------------------------------------

\end{document}